\documentclass[11pt, a4paper]{article}

\usepackage[all]{xy}

\usepackage[latin1]{inputenc}
\usepackage{amsfonts}
\usepackage{amsmath}
\usepackage{amssymb}
\usepackage{amscd}
\usepackage{amsthm}
\usepackage{indentfirst}

\usepackage{epsfig}

\usepackage{color}

\usepackage{cleveref}

\newtheorem{thm}{Theorem}[section]
\newtheorem{lemma}[thm]{Lemma}

\theoremstyle{definition}

\newtheorem{rem}[thm]{Remark}

\newtheorem{question}[thm]{Question}

\newtheorem*{acknow}{Acknowledgements}

\newtheorem*{prf}{Proof}

\newcommand{\R}{{\mathbb{R}}}

\newcommand{\Z}{{\mathbb{Z}}}
\newcommand{\N}{{\mathbb{N}}}

\newcommand{\cC}{{\mathcal{C}}}

\newcommand{\cQ}{{\mathcal{Q}}}

\newcommand{\om}{{\omega}}
\newcommand{\epsi}{{\varepsilon}}
\newcommand{\eps}{\varepsilon}

\newcommand{\fc}{{:\ }}

\newcommand{\vph}{\varphi}
\newcommand{\ol}{\overline}
\newcommand{\wt}{\widetilde}

\newcommand{\tb}{\textbf}

\DeclareMathOperator{\im}{im}

\DeclareMathOperator{\area}{area}
\DeclareMathOperator{\Vol}{Vol}

\DeclareMathOperator{\const}{const}

\DeclareMathOperator{\supp}{supp}

\DeclareMathOperator{\osc}{osc}

\DeclareMathOperator{\diam}{diam}

\newcommand{\til}{\tilde}

\newcommand{\setm}{\setminus}

\title{Rigidity of the $L^p$-norm of the Poisson bracket on surfaces}
\author{
	Karina Samvelyan\footnote{School of Mathematical Sciences, Faculty of Exact Sciences, Tel Aviv University, \texttt{karina.samvelyan@gmail.com}.}\ \
	and Frol Zapolsky\footnote{Department of Mathematics, Faculty of Natural Sciences, University of Haifa, \texttt{frol.zapolsky@gmail.com}.}}

\begin{document}

\renewcommand{\labelenumi}{(\roman{enumi})}

\maketitle

\begin{abstract}
For a symplectic manifold $(M,\omega)$ let $\{\cdot,\cdot\}$ be the corresponding Poisson bracket. In this note we prove that the functional
$$(F,G) \mapsto \|\{F,G\}\|_{L^p(M)}$$
is lower-semicontinuous with respect to the $C^0$-norm on $C^\infty_c(M)$ when $\dim M = 2$ and $p < \infty$, extending previous rigidity results for $p = \infty$ in arbitrary dimension.
\end{abstract}

\section{Introduction and main result}

One of the fascinating manifestations of rigidity in symplectic topology is the unexpected robust behavior of the Poisson bracket with respect to the $C^0$-norm on the space of smooth functions, discovered by Cardin--Viterbo \cite{cardin-viterbo}.
To state their seminal result, let $(M,\omega)$ be a symplectic manifold without boundary, and let us endow the space $C^\infty_c(M)$ of smooth compactly supported functions on $M$ with the topology induced by the supremum norm $\|\cdot\|_{C^0}$. We write $\xrightarrow{C^0}$ to indicate convergence with respect to this topology.

The Poisson bracket of $F,G \in C^\infty(M)$ is the function
$$\{F,G\} = -\omega(X_F,X_G) = dF(X_G)\,,$$
where for $H \in C^\infty(M)$ its Hamiltonian vector field $X_H$ is defined by $\omega(X_H,\cdot) = -dH$.
\begin{thm}[Cardin--Viterbo \cite{cardin-viterbo}]Let $N$ be $\R^n$ or a closed manifold,\footnote{The proof uses generating functions for Lagrangians in $T^*N$, therefore it is plausible that it extends to more general $N$, however this formulation suffices to illustrate the main point.} and assume that $M = T^*N$ and $\omega$ is the canonical symplectic form. Let $F,G \in C^\infty_c(M)$ be such that $\{F,G\} \neq 0$. Then
$$\liminf_{\ol F \xrightarrow{C^0},\,\ol G\xrightarrow{C^0} F}\|\{\ol F, \ol G\}\|_{C^0} > 0\,.$$
\end{thm}
\noindent This means that if two functions do not Poisson commute, it is impossible to approximate them, in the $C^0$ sense, by Poisson commuting, or even asymptotically commuting, functions. This behavior is surprising because the Poisson bracket is defined in terms of the first derivatives of the functions and thus \emph{a priori} it is unknown how it changes under $C^0$ perturbations. For surfaces, a stronger form of this statement was proved in \cite{Zapolsky_QS_PB_surfaces}: \footnote{See also \cite{Entov_Polterovich_Zapolsky_qms_Poisson_bracket} for intermediate quantitative results in arbitrary dimension using symplectic quasi-states.}

\begin{thm}Assume $\dim M = 2$. Then for $F,G \in C^\infty_c(M)$ the functional $\|\{\cdot,\cdot\}\|_{C^0}$ is lower-semicontinuous with respect to $C^0$-norm, meaning
$$\liminf_{\ol F \xrightarrow{C^0},\, \ol G \xrightarrow{C^0} G} \|\{\ol F, \ol G\}\|_{C^0} = \|\{F,G\}\|_{C^0}\,.$$
\end{thm}

\noindent This result was proved using methods of classical analysis in dimension two. It was later generalized, using methods of ``hard'' symplectic topology, including the Hofer metric and the energy-capacity inequality, to arbitrary dimension:
\begin{thm}[\cite{entov-polterovich2}, \cite{buhovsky}]For $M$ of arbitrary dimension and any $F,G \in C^\infty_c(M)$ we have
$$\liminf_{\ol F \xrightarrow{C^0} F,\,\ol G \xrightarrow{C^0} G} \|\{\ol F,\ol G\}\|_{C^0} = \|\{F,G\}\|_{C^0}\,.$$
\end{thm}

Our main result in this note is the following rigidity phenomenon in dimension two, proved using a refinement of the technique from \cite{Zapolsky_QS_PB_surfaces}:
\begin{thm}\label{thm:main_result}Assume $\dim M = 2$ and $p \in [1,\infty)$. Then for $F,G \in C^\infty_c(M)$ we have
$$\liminf_{\ol F \xrightarrow{C^0} F,\,\ol G \xrightarrow{C^0} G} \|\{\ol F,\ol G\}\|_{L^p(M)} = \|\{F,G\}\|_{L^p(M)}\,.$$
\end{thm}
\noindent Here and in the rest of the note we denote by $\|H\|_{L^p(X)}$ the $L^p$-norm, with respect to the measure induced by $\omega$, of a function $H$ defined on a measurable subset $X \subset M$.

Whether this behavior persists in higher dimension is currently unknown. Therefore we ask the following question.
\begin{question}
  Is the functional $(F,G) \mapsto \|\{F,G\}\|_{L^p(M)}$ lower semi-conti\-nuous for $M$ of arbitrary dimension and finite $p$?
\end{question}

\noindent In view of the results in \cite{buhovsky}, it is also natural to ask the following.
\begin{question}What is the modulus of semi-continuity of the functional $(F,G) \mapsto \|\{F,G\}\|_{L^p(M)}$? Is there a constant $\kappa > 0$ such that
$$\inf_{\|\ol F - F\|_{C^0},\,\|\ol G - G\|_{C^0} \leq \delta} \|\{\ol F,\ol G\}\|_{L^p(M)} \geq \|\{F,G\}\|_{L^p(M)} - \const(F,G)\cdot\delta^\kappa\,?$$
\end{question}

The rigid behavior with respect to the $C^0$-norm should be contrasted with the following result.
\begin{thm}[\cite{Samvelyan_msc_thesis}]\label{thm:PB_Lqp_flexibility}
  Let $M$ have arbitrary dimension $2n$, let $q \in [1,\infty)$, and let $F,G \in C^\infty_c(M)$. Then for any $\eps > 0$ and a compact submanifold with boundary $C \subset M$ of dimension $2n$, whose interior contains $\supp F \cup \supp G$, there exist $\wt F,\wt G \in C^\infty_c(M)$ supported in $C$, such that
    $$\|\wt F - F\|_{C^0} < \eps\,,\quad \|\wt G - G\|_{L^q(M)} < \eps^{1/q}\,,\quad \text{and} \quad \{\wt F,\wt G\} \equiv 0\,.$$
  In particular, for any $p \in [1,\infty]$,
    $$
		\liminf_{\wt F \xrightarrow{C^0} F,\,\wt G \xrightarrow{L^q} G} \|\{\wt F,\wt G\}\|_{L^p(M)} = \liminf_{\wt F \xrightarrow{L^q} F,\,\wt G \xrightarrow{L^q} G} \|\{\wt F,\wt G\}\|_{L^p(M)} = 0 \,. 
	$$
\end{thm}
\noindent Here we use the fact that $\|H\|_{L^q(M)} \leq \big(\int_C\omega^n\big)^{1/q}\cdot \|H\|_{C^0}$ for $H \in C^\infty(M)$ with $\supp H \subset C$. This means that the $L^p$-norm of the Poisson bracket becomes flexible if we take the $L^q$-topology on $C^\infty_c(M)$ for finite $q$.

For the sake of completeness, we provide a proof of the theorem in the next section.

\begin{rem}
Note that for continuous functions the $L^\infty$-norm and the $C^0$-norm coincide.
\end{rem}

\begin{acknow}We wish to thank Lev Buhovsky and Leonid Poltero\-vich for reading a preliminary version of the paper and making useful comments, and for their interest.
KS is partially supported by the Israel Science Foundation grant number 178/13, and by the European Research Council Advanced grant number 338809.
FZ is partially supported by grant number 1281 from the GIF, the German--Israeli Foundation for Scientific Research and Development, and by grant number 1825/14 from the Israel Science Foundation.
\end{acknow}

\section{Proofs}

\begin{prf}[of Theorem \ref{thm:main_result}]
Let us give an overview of the proof before passing to the details. The actual logical order of the proof is somewhat different from this summary.

We define the map
$$\Phi \fc M \to \R^2\quad \text{by} \quad \Phi(z) = (F(z),G(z))\,.$$
The main point is that since $\dim M = 2$, the Poisson bracket $\{F,G\}$ is related to $\Phi$ via
$$\Phi^*(dx \wedge dy) = dF \wedge dG = -\{F,G\}\omega\,,$$
where $(x,y)$ are the coordinate functions on $\R^2$. We see that a point $z \in M$ is regular for $\Phi$ if and only if $\{F,G\}(z) \neq 0$. We let $U \subset \R^2$ be the set of regular values of $\Phi$ in $\im \Phi$.

Consider now the subset $K_n \subset U$ comprised of squares of size $\frac 1 n$ with vertices in the grid $\frac 1 n \Z \times \frac 1 n \Z$ with $n \in \N$ large so that $\|\{F,G\}\|_{L^p(\Phi^{-1}(K_n))}$ is close to $\|\{F,G\}\|_{L^p(M)}$. Next we \emph{subdivide each square in $K_n$ into squares of size} $\frac 1 {kn}$, $k \in \N$. For such a square $Q$, a connected component $Q' \subset \Phi^{-1}(Q)$, and $k$ large, the \emph{oscillation of $\{F,G\}$ over $Q'$ can be made arbitrarily small} for all such $Q'$, which allows us to relate the $L^1$- and the $L^p$-norms of $\{F,G\}$ over $Q'$. Then we use the \emph{lower semi-continuity of the $L^1$-norm} to pass to $\|\{\ol F,\ol G\}\|_{L^1(Q')}$. Finally the \emph{H\"older inequality} brings us back to $\|\{\ol F,\ol G\}\|_{L^p(Q')}$.

\begin{rem}
We wish to note here that the use of two scales, $\frac 1 n$ and $\frac 1 {kn}$, seems to stem from convenience rather than being a reflection of something deeper. We must simultaneously approximate the $L^p$-norm of $\{F,G\}$ and control its oscillation, and this double subdivision is a way to do it.
\end{rem}

We now give the details of the proof.
\begin{rem} \label{rem:Phi_covering}
Since $M$ is assumed to have no boundary, and $\{F,G\}$ has compact support, $\Phi$ is a covering map over $U$. In particular, the lifting property of a covering implies that if $Y \subset U$ is a path-connected simply connected subset, then $\Phi|_{\Phi^{-1}(Y)} \fc \Phi^{-1}(Y) \to Y$ is a trivial covering, that is $\Phi^{-1}(Y)$ is a disjoint union of path components, each one projected homeomorphically onto $Y$ by $\Phi$. If $Y$ is in addition a submanifold with corners, then, since $\Phi$ is smooth, these components are themselves submanifolds with corners, projected in fact diffeomorphically onto $Y$.
\end{rem}

For $n \in \N$ let $K_n \subset \R^2$ be the union of squares of the form $[\frac{i}{n},\frac{i+1}{n}] \times [\frac{j}{n},\frac{j+1}{n}]$, where $i,j\in\Z$, contained in $U$. For $k \in \N$ consider a square $Q = [\frac{i}{kn},\frac{i+1}{kn}] \times [\frac{j}{kn},\frac{j+1}{kn}]$, where $i,j \in \Z$, and assume it is contained in $K_n$. By Remark \ref{rem:Phi_covering}, $\Phi^{-1}(Q)$ is a disjoint union of connected components, each of which is mapped by $\Phi$ diffeomorphically onto $Q$. See fig. \ref{fig: construction_of_K_n_and_cQ_k}. Let $\cQ_{n,k}$ be the collection of all such connected components for all such $Q$. The next lemma states that the oscillation of $|\{F,G\}|^p$ over the sets in $\cQ_{n,k}$ can be made arbitrarily small as $k \to \infty$.

\begin{figure}[h!]
	\centering
	\includegraphics[scale=0.9]{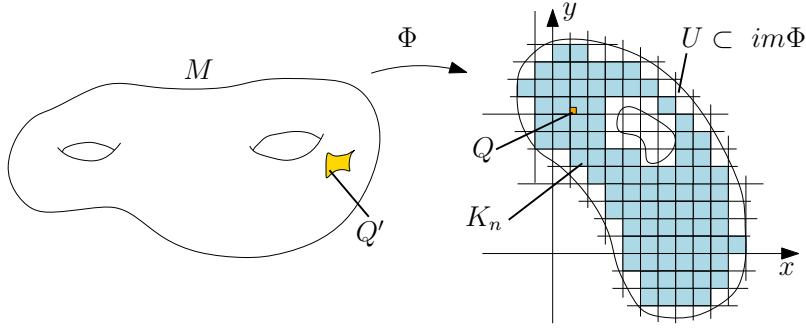}
	\caption{The set $K_n\subset U$ and an element $\Phi (Q') = Q$ in its subdivision.}
	\label{fig: construction_of_K_n_and_cQ_k}
\end{figure}

\begin{lemma}\label{lem:oscillation}
$\lim_{k \to \infty} \max_{Q' \in \cQ_{n,k}}\osc_{Q'}|\{F,G\}|^p = 0$.
\end{lemma}
Fix $\eps > 0$ and let $k \in \N$ be such that $\max_{Q' \in \cQ_{n,k}}\osc_{Q'}|\{F,G\}|^p \leq \eps$. Pick $Q' \in \cQ_{n,k}$, let $Q = \Phi(Q') \subset \R^2$, and let $i,j\in\Z$ be such that $Q = [\frac{i}{kn},\frac{i+1}{kn}] \times [\frac{j}{kn},\frac{j+1}{kn}]$.
For $\delta \in (0, \tfrac 1 {2kn})$ denote $Q_\delta = [\frac{i}{kn}+\delta,\frac{i+1}{kn} - \delta] \times [\frac{j}{kn} + \delta,\frac{j+1}{kn} - \delta]$. The following lemma is a quantitative local surjectivity result for $C^0$-perturbations of $\Phi$. Its proof is an almost verbatim repetition of the one of Lemma 3.1 in \cite{Zapolsky_QS_PB_surfaces} and is omitted.
\begin{lemma}\label{lem:local_surjectivity}
Let $\delta \in (0,\frac{1}{2kn})$ and let $\ol F,\ol G \in C^\infty_c(M)$ be such that
$$\|\ol F - F\|_{C^0} \leq \delta \,,\; \|\ol G - G\|_{C^0} \leq \delta\,.$$
Define
$$\ol \Phi \fc M \to \R^2 \quad \text{by} \quad \ol\Phi(z) = (\ol F(z),\ol G(z))\,.$$
Then we have
$$\ol \Phi(Q') \supset Q_\delta\,. \qed$$
\end{lemma}
Fix $\delta \in (0,\frac{1}{2kn})$ and $\ol F,\ol G \in C^\infty_c(M)$ with $\|\ol F - F\|_{C^0} \leq \delta$, $\|\ol G - G\|_{C^0} \leq \delta$. Let $q$ be such that $1/p + 1/q = 1$. The H\"older inequality allows us to relate the $L^p$- and the $L^1$-norms of $\{\ol F,\ol G\}$:
$$\|\{\ol F,\ol G\}\|_{L^p(Q')}^p \geq \|\{\ol F,\ol G\}\|_{L^1(Q')}^p \|1\|_{L^q(Q')}^{-p}\,.$$
Let us define the function
$$n_{\ol\Phi} \fc \R^2 \to \N \cup \{0,\infty\}\,,$$
where $n_{\ol\Phi}(u)$ is the number of preimages of $u$ by \emph{the restriction of} $\ol\Phi$ \emph{to} $Q'$. Note that by Lemma \ref{lem:local_surjectivity} $n_{\ol\Phi}(u) \geq 1$ for every $u \in Q_\delta$. The so-called area formula from geometric measure theory \cite[Theorem 3.2.3]{Federer_Geometric_meas_thry} implies in our case the following identity:
$$\int_{Q'} |d\ol F \wedge d\ol G| = \int_{\R^2}n_{\ol \Phi}\,dx \wedge dy\,,$$
where for a $2$-form $\beta$ on $M$ we let $|\beta|$ denote the corresponding density. \footnote{This can be thought of as the nonnegative measure induced by $\beta$; if $f \in C^\infty(M)$ is such that $\beta = f\omega$, then $\int |\beta| \equiv \int |f|\omega$.}

Next, we relate the $L^1$-norms of $\{\ol F,\ol G\}$ and $\{F,G\}$ over $Q'$:
\begin{align*}
\|\{\ol F, \ol G\}\|_{L^1(Q')} &= \int_{Q'}|d\ol F \wedge d\ol G| && \text{by the definition of }\{\cdot,\cdot\}\\
&= \int_{\R^2}n_{\ol\Phi} \,dx \wedge dy && \text{by the area formula}\\
&\geq \int_{Q_\delta} dx \wedge dy && \text{since }n_{\ol\Phi}|_{Q_\delta} \geq 1\,.
\end{align*}
The last integral is the area of $Q_\delta$, which equals
$$\big(\tfrac 1 {kn} - 2\delta\big)^2 = (1-2kn\delta)^2\area(Q) = (1-2kn\delta)^2 \int_Q dx \wedge dy\,.$$
We continue:
\begin{align*}
\|\{\ol F, \ol G\}\|_{L^1(Q')} &\geq (1-2kn\delta)^2 \int_Q dx \wedge dy \\
&= (1-2kn\delta)^2 \int_{\Phi(Q')} |dx \wedge dy|\\
&= (1-2kn\delta)^2 \int_{Q'} |\Phi^*(dx \wedge dy)| \\
&= (1-2kn\delta)^2 \int_{Q'}|dF \wedge dG|\\
&= (1-2kn\delta)^2 \|\{F,G\}\|_{L^1(Q')} \,,
\end{align*}
therefore
$$\|\{\ol F, \ol G\}\|_{L^1(Q')}^p \geq (1-2kn\delta)^{2p} \|\{F,G\}\|_{L^1(Q')}^p\,.$$
Now we relate the $L^1$- and the $L^p$-norms of $\{F,G\}$ over $Q'$. Namely, since $\osc_{Q'}|\{F,G\}|^p \leq \eps$, we have
\begin{multline*}
  \textstyle\|\{F,G\}\|_{L^1(Q')}^p \geq \big(\min_{Q'} |\{F,G\}| \int_{Q'} \omega \big)^p = \min_{Q'}|\{F,G\}|^p \big(\int_{Q'}\omega\big)^p \geq \\ \textstyle \geq (\max_{Q'}|\{F,G\}|^p - \eps) \big(\int_{Q'}\omega\big)^p\,,
\end{multline*}
therefore, since $\|1\|_{L^q(Q')}^{-p} = \big(\int_{Q'}\omega \big)^{-p/q}$:
$$\textstyle\|\{F,G\}\|_{L^1(Q')}^p \|1\|_{L^q(Q')}^{-p} \geq (\max_{Q'}|\{F,G\}|^p - \eps) \big(\int_{Q'}\omega\big)^{p-p/q}\,.$$
Since $p-p/q = 1$, we obtain
$$\textstyle \max_{Q'}|\{F,G\}|^p \big(\int_{Q'}\omega\big)^{p-p/q} = \max_{Q'}|\{F,G\}|^p \int_{Q'}\omega \geq \|\{F,G\}\|_{L^p(Q')}^p\,,$$
thus in total
$$\|\{F,G\}\|_{L^1(Q')}^p \|1\|_{L^q(Q')}^{-p} \geq \|\{F,G\}\|_{L^p(Q')}^p - \eps\int_{Q'}\omega\,.$$

Assembling all of the above, we obtain the main estimate
\begin{equation}\label{eqn:main_estimate}
  \textstyle\|\{\ol F,\ol G\}\|_{L^p(Q')}^p \geq (1-2kn\delta)^{2p}\big(\|\{F,G\}\|_{L^p(Q')}^p - \eps\int_{Q'}\omega\big)\,.
\end{equation}
Note that $\Phi^{-1}(K_n)$ is the essentially disjoint\footnote{A countable union of subsets is essentially disjoint if the intersection of every two subsets has measure zero.} union of the sets $Q' \in \cQ_{n,k}$, and that $\|\cdot \|_{L^p}^p$ is additive with respect to essentially disjoint unions. Thus we have
\begin{align*}
\|\{\ol F,\ol G\}\|_{L^p(M)}^p & \geq \|\{\ol F,\ol G\}\|_{L^p(\Phi^{-1}(K_n))}^p\\
&= \sum_{Q'\in\cQ_{n,k}}\|\{\ol F,\ol G\}\|_{L^p(Q')}^p \\
&\overset{*}{\geq} (1-2kn\delta)^{2p}\sum_{Q' \in \cQ_{n,k}}\textstyle \big(\|\{F,G\}\|_{L^p(Q')}^p - \eps\int_{Q'}\omega\big) \\
&= \textstyle (1-2kn\delta)^{2p} \big ( \|\{F,G\}\|_{L^p(\Phi^{-1}(K_n))}^p - \eps\int_{\Phi^{-1}(K_n)}\omega\big) \\
&\geq (1-2kn\delta)^{2p} \big(\|\{F,G\}\|_{L^p(\Phi^{-1}(K_n))}^p - \eps\cdot \area \supp \{F,G\}\big)\,,
\end{align*}
where for $\overset{*}{\geq}$ we used the main estimate \eqref{eqn:main_estimate}, and in the last inequality we used $\Phi^{-1}(K_n) \subset \supp \{F,G\}$. Taking $\delta \to 0$, we see that
$$\liminf_{\ol F \xrightarrow{C^0} F,\,\ol G \xrightarrow{C^0}G} \|\{\ol F,\ol G\}\|_{L^p(M)}^p \geq \|\{F,G\}\|_{L^p(\Phi^{-1}(K_n))}^p - \eps \cdot \area \supp \{F,G\}\,,$$
and since $\eps$ was arbitrary, we have
$$\liminf_{\ol F \xrightarrow{C^0} F,\,\ol G \xrightarrow{C^0}G} \|\{\ol F,\ol G\}\|_{L^p(M)}^p \geq \|\{F,G\}\|_{L^p(\Phi^{-1}(K_n))}^p\,.$$
It remains to invoke the following lemma, which says that the $L^p$-norm of $\{F,G\}$ can be approximated by looking at the sets $\Phi^{-1}(K_n)$.
\begin{lemma}\label{lem:Lp_norm_approxd_Kn}
$\sup_{n \in \N} \|\{F,G\}\|^p_{L^p(\Phi^{-1}(K_n))} = \|\{F,G\}\|_{L^p(M)}^p$.
\end{lemma}
\noindent The proof is thus finished, assuming Lemmas \ref{lem:oscillation}, \ref{lem:Lp_norm_approxd_Kn}. \qed
\end{prf}

It remains to prove the lemmas. We keep the notations introduced during the proof of Theorem \ref{thm:main_result}.

\begin{prf}[of Lemma \ref{lem:oscillation}]
Let $C \subset K_n$ be a square entering the definition of $K_n$. By Remark \ref{rem:Phi_covering}, $\Phi^{-1}(C)$ is a disjoint union of a finite number of components, each projecting diffeomorphically onto $C$ by $\Phi$. Let $\cC$ be the collection of all such connected components for all the squares $C \subset K_n$. Note that $\cC$ is finite. Pick $C' \in \cC$, let $C = \Phi(C')$, and let $P_{C'} \fc C \to \R$ be the function $|\{F,G\}|^p\circ (\Phi|_{C'})^{-1}$. Since $\Phi|_{C'} \fc C' \to C$ is a diffeomorphism, we have for any $Z \subset C'$:
$$\osc_{Z}|\{F,G\}|^p = \osc_{\Phi(Z)}P_{C'}\,.$$
It then follows that it is enough to prove the following for every $C' \in \cC$:
$$\lim_{k \to \infty}\;\max_{Q' \in \cQ_{n,k}\,,Q' \subset C'}\;\osc_{\Phi(Q')}P_{C'} = 0\,.$$
This follows from the fact that $P_{C'}$ is a smooth function, in particular it has bounded derivatives, and therefore its oscillation over $\Phi(Q')$ is bounded by a constant times the diameter of $\Phi(Q')$ which is $\frac {\sqrt 2}{kn}$. \qed
\end{prf}

\begin{prf}[of Lemma \ref{lem:Lp_norm_approxd_Kn}] \footnote{We thank Lev Buhovsky for a suggestion that lead to a simplification of the proof of the lemma.}
Let $X = \supp F \cap \supp G$, $V = \Phi^{-1}(U)$, and $Z = X - V$, which is the subset of $X$ consisting of points lying over singular values of $\Phi$. Let $S, R \subset M$ be the sets of critical and regular points of $\Phi$, respectively. We have the disjoint union \footnote{Note that there may be regular points of $\Phi$ which are mapped to singular values.}
$$Z = (Z \cap S) \cup (Z \cap R)\,.$$
At the beginning of the proof of Theorem \ref{thm:main_result} we noted that $z \in S$ if and only if $\{F,G\}(z) = 0$, therefore
$$\int_{Z \cap S}|\{F,G\}|^p\omega = 0\,.$$
We claim that $Z \cap R$ has measure zero. Indeed, $Z \cap R = (\Phi|_R)^{-1}(\im \Phi - U)$, and the claim follows from the fact that $R$ is an open subset of $M$, therefore a submanifold, $\Phi|_R$ is a local diffeomorphism, the fact that $\im \Phi - U$ has measure zero by Sard's theorem, and the following lemma.
\begin{lemma}Let $N,P$ be manifolds, let $f \fc N \to P$ be a local diffeomorphism, and let $Y \subset P$ a subset of measure zero. Then $f^{-1}(Y)$ has measure zero.
\end{lemma}
\begin{prf}
Since our manifolds are paracompact, they are second countable, and in particular $N$ can be covered with countably many charts, such that on each one of them $f$ is a diffeomorphism onto its image. Since diffeomorphisms preserve the property of having measure zero, it follows that $f^{-1}(Y)$ is covered by countably many measure zero sets, and thus it is itself such. \qed
\end{prf}
\noindent This implies
$$\int_{Z \cap R}|\{F,G\}|^p \omega = 0\,,$$
and therefore we have
$$\int_M |\{F,G\}|^p \omega = \int_X |\{F,G\}|^p \omega = \int_V |\{F,G\}|^p \omega\,.$$
From the regularity of the measure $|\{F,G\}|^p\omega$ we obtain
$$\int_V |\{F,G\}|^p \omega = \sup_{K \subset V\text{ compact}}\int_K |\{F,G\}|^p \omega\,.$$
It is therefore enough to show that for any compact $K \subset V$ there is $n \in \N$ such that $\Phi(K) \subset K_n$. This follows from the fact that $\Phi(K)$ is compact and contained in $U$, therefore $d(\Phi(K),\R^2 - U) > 0$ and
$$\lim_{n \to \infty} d(K_n,\R^2 - U) = 0\,,$$
where $d$ is the Euclidean distance between subsets of $\R^2$. This limit is indeed zero since $K_n$ contains all the points at a distance at least $\sqrt 2/n$ from $\R^2 - U$. \qed
\end{prf}

We now prove the flexibility result, Theorem \ref{thm:PB_Lqp_flexibility}.
\begin{prf}[of Theorem \ref{thm:PB_Lqp_flexibility}]
Let $C \subset M$ be as in the formulation of the theorem and fix $\eps > 0$. We need to construct a Poisson commuting pair $\wt F,\wt G \in C^\infty_c(M)$ supported in $C$ and satisfying
$$\|\wt F  - F\|_{C^0} < \eps\,, \quad \|\wt G - G\|_{L^p}^p < \eps\,.$$

Fix a Riemannian metric $d$ on $M$. By a simplex in $M$ we mean the image of an embedding $\Delta \to M$, where $\Delta$ is a closed simplex in $\R ^{2n}$. A triangulation of $C$ is a representation of $C$ as a union of such simplices, where every two simplices intersect only in a common face (which is a simplex of lower dimension). A construction described in \cite{Cairns_triangulation_1961} produces a finite such triangulation; moreover given $\delta > 0$, every simplex in this triangulation may be assumed to have diameter $< \delta$ with respect to $d$.

Since $C$ is compact, $F$ is uniformly continuous on it, that is there exists $\delta$ such that if $x,y \in C$ satisfy $d(x,y) < \delta$, then $|F(x) - F(y)| < \delta$. Fix such $\delta$ and take a triangulation of $C$ with all the simplices having diameter $< \delta$.

For every simplex $Q$ from the triangulation we fix open subsets with smooth boundary $Q_3 \Subset Q_2 \Subset Q_1 \Subset Q$, satisfying
$$\Vol (Q \setm Q_3) < \eps\cdot \frac{\Vol(Q)}{\|G\|_{C^0}^q\cdot\Vol(C)}\,.$$
Here $A \Subset B$ means that the closure of $A$ is contained in the interior of $B$, and $\Vol$ is the volume with respect to $\om^n$. The condition on the volumes is essential for constructing a suitable $\wt{G}$.

\tb{Construction of $\wt{F}$.} Consider a simplex $Q$ with open subsets $Q_2 \Subset Q_1 \Subset Q$ as above. We take an auxiliary smooth function $\vph \fc Q \to [0,1]$ such that $\vph |_{Q_2} \equiv 0$ and $\vph |_{Q \setm Q_1} \equiv 1$. Fix a point $x_0 \in Q_2$. Define $\wt{F}$ on $Q$ to be
$$\wt{F} (x) = \vph (x) F(x) + (1 - \vph (x)) F(x_0) \,.$$
We see that on $Q_2$ we have $\wt{F} \equiv F(x_0)$, while outside $Q_1$ we have $\wt{F} \equiv F$. See \cref{fig: construction_F_G_tildes_PB_flexibility}. Next, define $\wt F$ on $C$ by gluing all these partially defined functions. Note that the resulting function is well-defined and smooth. Moreover, since $F$ vanishes near $\partial C$, it is also true for $\wt F$. Therefore we can extend $\wt F$ by zero to a smooth function on $M$ with support in $C$. For any $x \in Q$ we have
\begin{align*}
|\wt{F} (x) - F(x)| & = |\vph(x)F(x) + (1-\vph(x))F(x_0) - F(x)| =\\
& =  \underbrace{|1-\vph(x)|}_{\leq 1} \cdot \underbrace{|F(x)-F(x_0)|}_{< \epsi}
< \eps \,,
\end{align*}
where the last inequality holds since $\diam (Q) < \delta$. Since $Q$ is arbitrary, we obtain $\| \wt F- F \|_{C^0} < \eps$.

\tb{Construction of $\wt{G}$.} Consider again a simplex $Q$ from our triangulation with the subsets $Q_3 \Subset Q_2 \Subset Q_1 \Subset Q$, such that
$$\Vol (Q \setm Q_3) \leq \eps\cdot \frac{\Vol(Q)}{\|G\|_{C^0}^q\cdot\Vol(C)}\,.$$
Take a smooth function $\psi \fc Q \to [0,1]$ satisfying $\psi |_{Q_3} \equiv 1$, $\psi |_{Q \setm Q_2} \equiv 0$, and define $\wt G \fc Q \to \R$ by $\wt{G} = \psi G$. We have $\wt{G} \equiv G$ on $Q_3$ and $\til{G} |_{Q \setm Q_2} \equiv 0$. 
Take $\wt G$ to be the function on $M$ defined in this way on every simplex $Q$, and extended by zero to $M \setm C$. This again is a well-defined smooth function with support in $C$.

\begin{figure} [h!]
	\centering
	\includegraphics[scale=0.9]{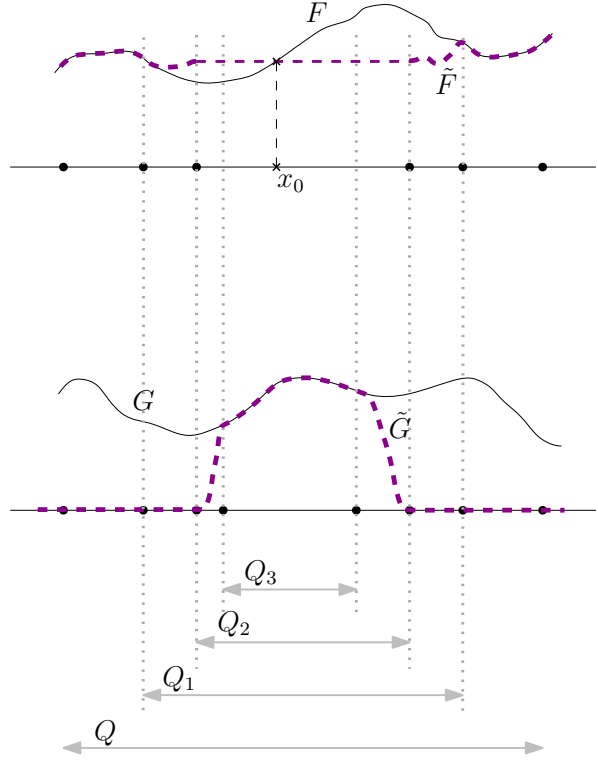}
	\caption{Producing $\wt F$ and $\wt G$ (the dashed lines).}
	\label{fig: construction_F_G_tildes_PB_flexibility}
\end{figure}

On a single simplex $Q$ we have
\begin{align*}
\int_Q | \wt{G} - G |^q \om^n
&= \int_{Q\setm Q_3} | \wt{G} - G |^q \om^n\\
&= \int_{Q\setm Q_3} (1-\psi)^q|G|^q \om^n\\
& \leq \int_{Q \setm Q_3} |G|^q \om^n\\
&\leq \|G\|_{C^0}^q\cdot \Vol(Q\setm Q_3)\\
& < \|G\|_{C^0}^q \cdot \eps \cdot \frac {\Vol(Q)}{\|G\|_{C^0}^q\Vol(C)} = \eps \cdot \frac{\Vol(Q)}{\Vol(C)}\,.
\end{align*}
\noindent Therefore on the whole of $M$ we get the bound
$$ \| \wt{G} - G\|_{L^q(M)}^q = \int_M | \wt{G} - G |^q \om^n < \eps \,.$$

It remains to note that for every simplex $Q$ in our triangulation and the associated subsets $Q_3 \Subset Q_2 \Subset Q_1 \Subset Q$, $\wt F$ is constant on $Q_2$, while $\wt G \equiv 0$ outside $Q_2$, meaning $\{\wt F,\wt G\} \equiv 0$ as claimed. \qed
\end{prf}

\bibliography{biblio}
\bibliographystyle{alpha}

\end{document}